\documentclass[12pt, en]{article}
\usepackage[margin=1.5cm]{geometry}                		% See geometry.pdf to learn the layout options. There are lots.
\usepackage[table]{xcolor}

\usepackage{authblk}
\usepackage{abstract}
\usepackage{booktabs}
\usepackage{enumitem}
\usepackage{natbib}
\usepackage{amssymb, amsthm, bbold, amsmath} 
\usepackage{algorithm}
\usepackage[noend]{algpseudocode}
\usepackage{graphicx}	
\usepackage[font=large]{subcaption}
\usepackage[displaymath, mathlines]{lineno}
\usepackage{xcolor}
\usepackage{hyperref}

\newcommand{\bx}{\bold{x}}
\newcommand{\bdx}{\delta{\bold{x}}}
\newcommand{\bL}{\bold{L}}
\newcommand{\bB}{\bold{B}}
\newcommand{\bR}{\bold{R}}
\newcommand{\bQ}{\bold{Q}}
\newcommand{\by}{\bold{y}}
\newcommand{\bH}{\bold{H}}
\newcommand{\bM}{\bold{M}}
\newcommand{\bI}{\bold{I}}

\newcommand{\bz}{\bold{z}}
\newcommand{\bC}{\bold{C}}

\title{On time-parallel preconditioning for the state formulation of incremental weak constraint 4D-Var}

\author[1]{\normalsize Ieva Dau\v{z}ickait\.{e}*}

\author[1,2]{\normalsize Amos S. Lawless}

\author[1,3]{\normalsize Jennifer A. Scott}

\author[1,4]{\normalsize Peter Jan van Leeuwen}

\affil[1]{\footnotesize School of Mathematical, Physical and Computational Sciences, University of Reading,UK}

\affil[2]{ \footnotesize National Centre for Earth Observation, Reading, UK}

\affil[3]{\footnotesize Scientific Computing Department, STFC Rutherford Appleton Laboratory, UK}

\affil[4]{\footnotesize Department of Atmospheric Science, Colorado State University, USA}

\date{}

\begin{document}

\maketitle

%\nolinenumbers
\begin{abstract}
Using a high degree of parallelism is essential to perform data assimilation efficiently. The state formulation of the incremental weak constraint four-dimensional variational data assimilation method allows parallel calculations in the time dimension. In this approach, the solution is approximated by minimising a series of quadratic cost functions using the conjugate gradient method. To use this method in practice, effective preconditioning strategies that maintain the potential for parallel calculations are needed. We examine approximations to the control variable transform (CVT) technique when the latter is beneficial. The new strategy employs a randomised singular value decomposition and retains the potential for parallelism in the time domain. Numerical results for the Lorenz 96 model show that this approach accelerates the minimisation in the first few iterations, with better results when CVT performs well. 
\vspace{\baselineskip}


\end{abstract}

\let\thefootnote\relax\footnotetext{*Correspondence: Ieva Dau\v{z}ickait\.{e}, Department of Mathematics and Statistics, Pepper Lane, Whiteknights, Reading RG6 6AX, UK. Email: i.dauzickaite@pgr.reading.ac.uk}

%\linenumbers

\section{Introduction}
The ever increasing resolution of weather models enhances the importance of parallelisation in data assimilation. Higher potential for parallel computations can be achieved by using suitable data assimilation methods. The state formulation of the weak constraint 4D-Var method, which allows for the model error, is such a method. In its incremental version, a series of quadratic cost functions is minimised via solving a series of linear systems containing the Hessian of the linearised cost function. These are solved with the conjugate gradient (CG) method (e.g. \cite{SaadBook}), where the most computationally expensive part is integrating the tangent linear model and its adjoint. It has been shown that these calculations can be parallelised in the time dimension (\cite{Fisher17}).

However, CG needs preconditioning for fast convergence. Efficient preconditioning for the state formulation of incremental weak constraint 4D-Var, which also preserves the potential for parallel in time calculations, is still an open question. By analogy with the standard preconditioning technique (also known as a control variable transform or first level preconditioning) used in strong constraint 4D-Var, \cite{Fisher17} suggested using approximations of the tangent linear model. Their search for a suitable approximation was unsuccessful. Our investigation in this paper reveals that preconditioning using the exact tangent linear model can be detrimental to the minimisation in some cases. We focus on approximations in the case when using the exact tangent linear model works well.

In the light of the growing popularity of randomised methods and examples of their use in data assimilation (\cite{Bousserez2020, Dauzickaite2021}), we propose using a randomised singular value decomposition (RSVD) (\cite{Halko11}) to approximate the tangent linear model. RSVD is a block method that is easy to parallelise in the sense that it requires calculating matrix products with blocks of vectors. Because \cite{Lawless2008} showed that it is important to take into account the information on the background errors when using model reduction techniques in data assimilation, we also examine an approach where we approximate the tangent linear model in interaction with the background and model error covariance matrices.

We formulate the incremental weak constraint 4D-Var problem and discuss its preconditioning in Section~\ref{sec:incremental_wc4D-Var}. Our idea for randomised preconditioning is presented in Section~\ref{sec:randomised_prec}. Numerical experiments exploring preconditioning using the exact tangent linear model and its low rank approximation obtained using RSVD are presented in Section~\ref{sec:numerics} and we summarize our findings and suggest future directions in Section~\ref{sec:conclusions}.

\section{Incremental weak constraint 4D-Var}\label{sec:incremental_wc4D-Var}

In data assimilation, the prior estimate of a model trajectory is combined with observations over a time window to obtain an improved estimate of the state (analysis) $\bx^a_0, \bx^a_1, \dots, \bx^a_N$ at times $t_0, t_1, \dots, t_N$. The prior estimate of the state at $t_0$ is called the background and is denoted by $\bx^b \in \mathbb{R}^n$ and the observations at time $t_i$ are denoted by $\by_i \in \mathbb{R}^{p_i}$. The state variables $\bx_i$ are mapped to the observation space using an observation operator $\mathcal{H}_i$. The nonlinear dynamical model $\mathcal{M}_i$ describes the state evolution from time $t_i$ to $t_{i+1}$. It is assumed that the background, observations, and model have Gaussian errors with zero mean and covariance matrices $\bB \in \mathbb{R}^{n \times n}$, $\bR_i \in \mathbb{R}^{p_i \times p_i}$, and $\bQ_i \in \mathbb{R}^{n \times n}$, respectively. We assume that the observation and model errors are uncorrelated in time.

In the state formulation of weak constraint 4D-Var, the analysis is the minimiser of the nonlinear cost function
\begin{equation}\label{eq:nonlinear-wc-4D-var}
J(\bx_0, \bx_1, \dots, \bx_N)  =  \frac{1}{2} \| \bx_0 - \bx^b\|^2 _{\bB^{-1}}+ \frac{1}{2} \sum_{i=0}^{N} \| \by_i - \mathcal{H}_i (\bx_i)\|^2_{\bR_i^{-1}}  
 + \frac{1}{2} \sum_{i=0}^{N-1} \| \bx_{i+1} - \mathcal{M}_i (\bx_i) \|^2_{\bQ_{i+1}^{-1} },
\end{equation}
where $\|\bold{a}\|^2_\bold{A} = \bold{a}^T \bold{A} \bold{a}$ (\cite{Tremolet06}).

The minimiser of \eqref{eq:nonlinear-wc-4D-var} can be approximated using an inexact Gauss-Newton algorithm (\cite{Gratton2007}). In this incremental approach, the $(j+1)$-st approximation \linebreak
$\bx^{(j+1)} = ( \bx_0^{(j+1) T}, \bx_1^{(j+1) T}, \dots, \bx_N^{(j+1) T})^T \in \mathbb{R}^{(N+1)n}$ of the state is
\begin{equation}
\bx^{(j+1)} = \bx^{(j)} + \delta \bx^{(j)},
\end{equation}
where the update is $\delta \bx^{(j)} = ( \delta \bx_0^{(j) T}, \delta \bx_1^{(j) T}, \dots, \delta \bx_N^{(j) T})^T \in \mathbb{R}^{(N+1)n}$. $\bM_i$ and $\bH_i$ are the model and observation operators linearised at $\bx_i$; they are known as the tangent linear model and tangent linear observation operator, respectively. %This approach uses the tangent linear model $\bM_i$ and observation operator $\bH_i$ linearised at time $t_i$. 
We define the following matrices (following \cite{Gratton18})\par
\begin{align}
 \bL^{(j)} = & \left( \begin{array}{ccccc}
\bI         &             &            &                    & \\
-\bM_0^{(j)} &     \bI       &            &                    &  \\
          &   -\bM_1^{(j)}&       \bI     &                    &  \\
         &             & \ddots & \ddots        &  \\
         &            &              & -\bM_{N-1}^{(j)}   & \bI  \\
\end{array} \right) \in \mathbb{R}^{(N+1)n \times (N+1)n}, \\
\bold{H}^{(j)} = &\  diag(\bH_0^{(j)}, \bH_1^{(j)}, \dots, \bH_N^{(j)} )\in \mathbb{R}^{p \times (N+1)n }, \\  
\bold{D}  = & \ diag(\bB,  \bQ_1, \dots,  \bQ_N) \in \mathbb{R}^{(N+1)n \times (N+1)n}, \\
\bold{R} = & \ diag(\bR_0, \bR_1, \dots, \bR_N) \in \mathbb{R}^{p \times p},
\end{align}
where $\bI \in \mathbb{R}^{n \times n}$ is the identity matrix, $diag(\cdot)$ denotes a block diagonal matrix and $p = \Sigma_{i=0}^N p_i$ is the total number of observations. We use the following notation for vectors\par
\begin{align}
\bold{b}^{(j)} = & \left( \begin{array}{c}
\bx_0^{(j)} - \bx^b\\
\mathcal{M}_0 (\bx_0^{(j)}) - \bx_1^{(j)} \\
\vdots \\
\mathcal{M}_{N-1} (\bx_{N-1}^{(j)}) - \bx_N^{(j)} 
\end{array} \right) \in \mathbb{R}^{(N+1)n}, \\ 
\bold{d}^{(j)} = & \left( \begin{array}{c}
\by_0  - \mathcal{H}_0 (\bx_0^{(j)})\\
\by_1  - \mathcal{H}_1 (\bx_1^{(j)}) \\
\vdots \\
\by_N  - \mathcal{H}_N (\bx_N^{(j)}) 
\end{array} \right) \in \mathbb{R}^p.
\end{align} 
The update $\delta \bx^{(j)} $ is the minimiser of
\begin{equation}\label{eq:incr_wc4d-var_state}
J^{\delta} (\delta \bx^{(j)}) = \frac{1}{2} || \bL^{(j)} \bdx^{(j)} - \bold{b}^{(j)} ||^2_{\bold{D}^{-1}} + \frac{1}{2} || \bold{H}^{(j)} \bdx^{(j)} - \bold{d}^{(j)} ||^2_{\bold{R}^{-1}}.
\end{equation}
Because \eqref{eq:incr_wc4d-var_state} is a quadratic cost function, $\delta \bx^{(j)} $ can be found by solving the following large linear systems with the Hessian $\bold{A}^{(j)} $ of $J^{\delta} (\delta \bx^{(j)})$ \par
\begin{align}
\bold{A}^{(j)}  \bdx^{(j)}  = & (\bL^T)^{(j)}  \bold{D}^{-1} \bold{b}^{(j)}   + (\bold{H}^T)^{(j)}   \bold{R}^{-1} \bold{d}^{(j)} , \label{eq:incr_wc4d-var_state_SPD_system}\\
\text{where } \bold{A}^{(j)}   = & (\bL^T)^{(j)}  \bold{D}^{-1} \bL^{(j)}  + (\bold{H}^T)^{(j)}   \bold{R}^{-1} \bold{H}^{(j)}. \label{eq:incr_wc4d-var_state_SPD_matrix}
\end{align}
It is assumed that $p\ll (N+1)n$, thus $(\bold{H}^T)^{(j)}   \bold{R}^{-1} \bold{H}^{(j)}$ is symmetric positive semi-definite. Because $(\bL^T)^{(j)}  \bold{D}^{-1} \bL^{(j)} $ is symmetric positive definite, $\bold{A}^{(j)}  \in \mathbb{R}^{(N+1)n \times (N+1)n }$ is symmetric positive definite. Hence the method of choice for solving \eqref{eq:incr_wc4d-var_state_SPD_system} is CG. Each iteration of CG requires one matrix-vector product with $\bold{A}^{(j)}$, which is expensive due to the tangent linear model and its adjoint in $\bL^{(j)}$ and $ (\bL^T)^{(j)}$, respectively. \cite{Fisher17} noted that the structure of $\bL^{(j)}$ allows the matrix-vector products with $\bold{A}^{(j)}$ to be parallelised in the time dimension, i.e. computation of $\bL^{(j)} \bold{z}$, where $\bold{z} \in \mathbb{R}^{(N+1)n}$, can be parallelised for the model linearised at different times. In the rest of this paper, the superscript ${(j)}$ is omitted for ease of notation.

In general, CG needs preconditioning for fast convergence. Efficient preconditioning maps the system to another system that can be solved faster and the solution of the original problem can be easily recovered from the solution of the preconditioned problem. Choosing a suitable preconditioner is highly problem dependent. Given the possibility of parallel computations in matrix-vector products with \eqref{eq:incr_wc4d-var_state_SPD_matrix}, the preconditioner should keep this potential.

\subsection{Preconditioning}

We consider an extension of the control variable transform, also known as a first level preconditioning, that is used in 3D-Var, where the model evolution is omitted, and in the strong constraint formulation of 4D-Var, where the model is assumed to have no error (see, e.g. \cite{Lorenc2000,Rawlins2007,Lawless2013}). The idea is to apply the preconditioner so that the first term of the preconditioned Hessian is equal to identity. Then the preconditioned Hessian is a sum of the identity matrix and a low rank symmetric positive semi-definite matrix with rank at most $p$. Its smallest eigenvalue is equal to one and it has at most $p$ eigenvalues that are larger than one. The latter can impair CG convergence if they are not well separated (for a more general discussion see e.g.\cite{Nocedal06,Liesen_Strakos}). 

Applying this kind of preconditioning to the state formulation of weak constraint 4D-Var requires preconditioning with $\bL^{-1} \bold{D}^{1/2}$, where
\begin{equation}\label{eq:Linv_full}
\bL^{-1}= \left( \begin{array}{ccccc}
\bI         &             &            &                    & \\
\bM_{0,0}    &     \bI       &            &                    &  \\
\bM_{0,1}   &   \bM_{1,1}   &       \bI     &                    &  \\
\vdots              &   \vdots        & \ddots & \ddots        &  \\
\bM_{0,N-1}  &   \bM_{1,N-1}         &       \cdots       & \bM_{N-1,N-1}   & \bI  \\
\end{array} \right)
\end{equation}
and $\bM_{i,j}=\bM_j \dots \bM_i$ denotes the linearised model integration from time $t_i$ to $t_{j+1}$. Matrix-vector products with $\bL^{-1}$ are sequential in the time dimension, i.e. $\bL^{-1} \bold{z} = (\bz_0^T, (\bM_0 \bz_0  + \bz_1)^T,  (\bM_1 (\bM_0 \bz_0  +   \bz_1)  + \bz_2)^T, \dots, (\bM_{N-1}( \bM_{N-2} \dots \bM_0 \bz_0 + \bM_{N-2}  \dots \bM_1 \bz_1 + \dots + \bz_{N-1})+ \bz_N)^T)^T$, where $\bold{z} =(\bz_0^T, \bz_1^T, \dots, \bz_N^T)^T$. \cite{Fisher17} suggested using an approximation $\widetilde{\bL}^{-1} $ of $\bL^{-1}$ in the preconditioner. Then the preconditioned system to be solved is \par
\begin{align}\label{eq:state_Linv_prec}
\bold{A}^{pr} \delta \widetilde{ \bx} & =\bold{D}^{1/2} \widetilde{\bL}^{-T} (\bL^T \bold{D}^{-1} \bold{b}  + \bold{H}^T  \bold{R}^{-1} \bold{d}), \\ 
\text{where} \qquad \bold{A}^{pr} &= \bold{D}^{1/2} \widetilde{\bL}^{-T} (\bL^T \bold{D}^{-1} \bL + \bold{H}^T  \bold{R}^{-1} \bold{H}) \widetilde{\bL}^{-1} \bold{D}^{1/2}, \label{eq:state_Aprec_mtx} \\
 \widetilde{\bL}^{-1} \bold{D}^{1/2} \delta \widetilde{ \bx} &= \bdx.
\end{align}
With an appropriate choice of $\widetilde{\bL}^{-1}$, $\bold{A}^{pr}$ is symmetric positive definite. $\widetilde{\bL}^{-1}$ should be chosen so that it can be applied in parallel. Fisher and G{\"u}rol could not find a suitable approximation that would guarantee good convergence. \cite{Gratton2018,Gratton18} discussed using $\widetilde{\bL}^{-1}$ where $\bM_i$ is set to zero or to the identity matrix in \eqref{eq:Linv_full}, which may be useful if the model state does not change significantly from one time step to the next. This may be unrealistic. We propose a new approximation strategy that avoids this assumption in the next section.

\section{Randomised preconditioning}\label{sec:randomised_prec}

Randomised methods for low rank matrix approximations have attracted a lot of interest in recent years because they require matrix products with blocks of vectors that can be easily parallelised and it has been shown that good approximations for matrices with rapidly decaying singular values can be obtained with high probability (e.g. \cite{Halko11}, \cite{Martinsson20}). These methods have been explored in data assimilation when designing solvers for strong constraint 4D-Var (\cite{Bousserez2020}) and preconditioning for the forcing formulation of the incremental weak constraint 4D-Var (\cite{Dauzickaite2021}). 

A low rank approximation of $\bL^{-1}$ cannot be used in \eqref{eq:state_Linv_prec}, because it would make \eqref{eq:state_Aprec_mtx} low rank and thus singular. Hence, we suggest exploiting the structure of $\bL^{-1}$ when generating the preconditioner. We write 
\begin{equation}\label{eq:Linv_I_plus_P}
\bL^{-1} = \bold{I} + \bold{P},
\end{equation}
where $\bold{P}$ is a strictly lower triangular matrix (with $0$ on the diagonal). We propose using a rank $k$ approximation $\widetilde{\bold{P}}=\bold{U} \bold{\Sigma} \bold{V}^T$, where $k$ is small compared to $(N+1)n$ and $\bold{U}, \bold{V} \in \mathbb{R}^{(N+1)n \times k}$, $\bold{\Sigma} \in \mathbb{R}^{k\times k}$ is a truncated singular value decomposition, i.e. $\bold{\Sigma}$ is diagonal with approximations to the $k$ largest singular values of $\bold{P}$ on the diagonal, and the columns of $\bold{U}$ and $\bold{V}$ are approximate left and right singular vectors, respectively. Then a non singular $\widetilde{\bL}^{-1}$ is %in \eqref{eq:state_Linv_prec} can be written as
\begin{equation}\label{eq:approxLinv}
\widetilde{\bL}^{-1} = \bold{I} + \widetilde{\bold{P}} = \bold{I} + \bold{U} \bold{\Sigma} \bold{V}^T.
\end{equation}
%\eqref{eq:approxLinv} is cheap to apply

A randomised singular value decomposition (RSVD) can be used to obtain $\widetilde{\bold{P}}$. RSVD is essentially one iteration of a classic subspace iteration method (\cite{Gu2015}). To increase the probability of success, the randomised methods work with larger subspaces than the required rank of the approximation. This is called oversampling. \cite{Halko11} indicate that setting the oversampling parameter $l$ to five or ten generally gives good results. We present the RSVD in Algorithm~\ref{alg:rsvd}. The entries of the Gaussian random matrix are independent standard normal random variables. Note that we remove the smallest $l$ computed singular values and the corresponding singular vectors. In this way the oversampling increases the cost of generating the preconditioner (in particular, the cost of the matrix-matrix products in steps~\ref{step:AG} and \ref{step:get_K} ), but not of its application. RSVD needs two expensive matrix-matrix products with $\bold{P}$ (steps~\ref{step:AG} and \ref{step:get_K} in Algorithm~\ref{alg:rsvd}), where $\bold{P}$ is multiplied with a matrix of size $(N+1)n \times (k+l)$. For efficiency, these can be parallelised. These matrix-matrix products consist of products with $\bM_{i,j}$. Hence, the cost of generating the preconditioner depends on the cost of integrating the tangent linear model over the assimilation window sequentially.

\cite{Lawless2008} showed that including the background error covariance matrix when using model reduction methods may lead to better results. Hence, we also explore using an approximation of 
\begin{equation}
\bold{S} = \bL^{-1} \bold{D}^{1/2} 
=\left( \begin{array}{ccccc}
\bB^{1/2}         &             &            &                    & \\
\bM_{0,0} \bB^{1/2}   &     \bQ_1^{1/2}       &            &                    &  \\
\bM_{0,1} \bB^{1/2} &   \bM_{1,1} \bQ_1^{1/2}    &       \bQ_2^{1/2}      &                    &  \\
\vdots              &   \vdots        & \ddots & \ddots        &  \\
\bM_{0,N-1} \bB^{1/2}  &   \bM_{1,N-1}  \bQ_1^{1/2}       &    \cdots          & \bM_{N-1,N-1} \bQ_{N-1}^{1/2}  & \bQ_N^{1/2}  \\
\end{array} \right). 
\end{equation}
As when approximating $\bL^{-1}$, we write 
\begin{equation}\label{eq:S_Dsqrt_plus_W}
\bold{S} = \bold{D}^{1/2} + \bold{W},
\end{equation}
where $\bold{W}$ is strictly lower triangular. Approximation $\widetilde{\bold{S}} = \bold{D}^{1/2} + \widetilde{\bold{W}}$ can be obtained by using RSVD to generate a low rank approximation of $\bold{W}$. The system to be solved is \eqref{eq:state_Linv_prec} with $\widetilde{\bL}^{-1} \bold{D}^{1/2}$ replaced by $\widetilde{\bold{S}}$.

\begin{algorithm}
\caption{Randomised singular value decomposition (RSVD)}\label{alg:rsvd}
\hspace*{\algorithmicindent} \textbf{Input:} matrix $\bold{A} \in \mathbb{R}^{s \times s}$,  target rank $k$, an oversampling parameter $l$ \\
\hspace*{\algorithmicindent} \textbf{Output:} orthogonal $\bold{U} \in \mathbb{R}^{s \times k}$ and $\bold{V} \in \mathbb{R}^{ s \times k}$  whose columns are approximations to left and right singular vectors of $\bold{A}$, respectively, and diagonal $\bold{\Sigma} \in \mathbb{R}^{ k \times k}$ with approximations to the largest singular values of $\bold{A}$
\begin{algorithmic}[1]
\State Form a Gaussian random matrix $\bold{G} \in \mathbb{R}^{s \times (k+l)}$
\State Form a sample matrix $\bold{Y}=\bold{A}\bold{G} \in \mathbb{R}^{s \times (k+l)}$ \label{step:AG}
\State Orthonormalize the columns of $\bold{Y}$ to obtain orthonormal $\bold{Z} \in \mathbb{R}^{s \times (k+l)}$
\State Form $\bold{K} = \bold{Z}^T \bold{A} \in \mathbb{R}^{ (k+l) \times s}$ \label{step:get_K}
\State Form SVD of $\bold{K}:\ \bold{K} = \hat{\bold{U}} \bold{\Sigma} \bold{V}^T$, where $\hat{\bold{U}}, \ \bold{\Sigma} \in \mathbb{R}^{ (k+l) \times (k+l)}, \ \bold{V} \in \mathbb{R}^{ s \times (k+l)}$
\State Remove last $l$ columns and rows of $\bold{\Sigma}$, so that $\bold{\Sigma} \in \mathbb{R}^{ k \times k}$
\State Remove last $l$ columns of $\hat{\bold{U}}$ and $\bold{V}$, so that $\hat{\bold{U}} \in \mathbb{R}^{ (k+l) \times k}, \bold{V} \in \mathbb{R}^{ s \times k}$
\State Form $\bold{U} = \bold{Z} \hat{\bold{U}} \in \mathbb{R}^{s \times k}$.
\end{algorithmic}
\end{algorithm} 

\section{Numerical results}\label{sec:numerics}
We test preconditioning using $\bL^{-1}$ and the approximations $\widetilde{\bold{L}}^{-1}$ and $\widetilde{\bold{S}}$ in \eqref{eq:state_Linv_prec} numerically. Preconditioning using the exact $\bL^{-1}$ is considered so that we understand when preconditioning using $\widetilde{\bL}^{-1}$ or $\widetilde{\bold{S}}$ may be effective but this is not regarded as a practical approach when parallelisation in the time dimension is desired. Identical twin experiments are performed. The background state $\bx^b$ is generated by adding random, Gaussian noise with covariance $\bB$ to $\bx_0^t$, where $\bx_i^t$ is the reference state at time $t_i$. We use direct observations that are obtained by adding random, Gaussian noise with covariance $\bR_i$ to $\mathcal{H}_i (\bx_i^t)$. 

The nonlinear Lorenz 96 model (\cite{Lorenz96}) is used, where the dynamics of $\bx_i = (X^1, \dots, X^n)^T$ are described by a set of $n$ coupled ODEs:
\begin{equation}\label{eq:lorenz96}
\frac{dX^j}{dt} = -X^{j-2} X^{j-1} + X^{j-1} X^{j+1} - X^j + F
\end{equation}
with conditions $X^{-1} = X^{n-1}, X^0 = X^n$ and $X^{n+1} = X^1$ and $F=8$. We use a fourth order Runge-Kutta scheme (\cite{Butcher87}). We consider the system with $n=100$ and $N=149$, so $\bold{A}^{pr}$ is a $15000 \times 15000$ matrix. The time step is set to $\Delta t = 2.5 \times 10^{-2}$ and the grid point distance is $\Delta X = 1/n$.

The covariance matrices are $\bB = 0.2^2 \bC_b$, $\bQ_i = 0.05^2 \bC_q$, where $\bC_b$ is a second-order autoregressive (SOAR) (\cite{Daley91}) matrix and $\bC_q$ a Laplacian (\cite{Johnson2005}) correlation matrix with length scales $2\Delta X$ and $0.75 \Delta X$, respectively. We consider $\bR_i=\sigma_o^2 \bI$ and vary $\sigma_o$.

The computations are performed with Matlab R2019b and the linear systems are solved with the Matlab preconditioned conjugate gradient (PCG) implementation \textit{pcg}.

\subsection{Preconditioning with exact $\bL^{-1}$}\label{sec:prec_exact_Linv}
\begin{figure}[h!]
\begin{subfigure}[b]{0.5\linewidth}
  \centering
 \includegraphics[width=\linewidth]{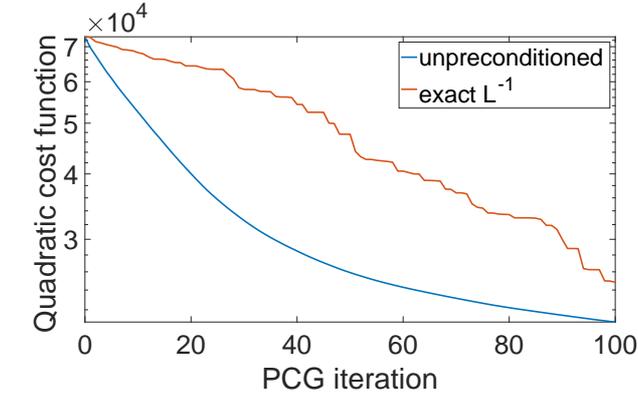}
 \caption{Case~\ref{case:reference}}
\end{subfigure}
\begin{subfigure}[b]{0.5\linewidth}
  \centering
 \includegraphics[width=\linewidth]{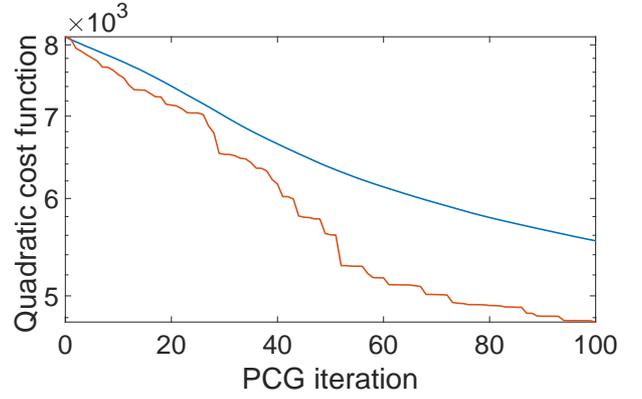}
  \caption{Case~\ref{case:large_sigmao}}
\end{subfigure}\\[1ex]
\begin{subfigure}[b]{0.5\linewidth}
  \centering
 \includegraphics[width=\linewidth]{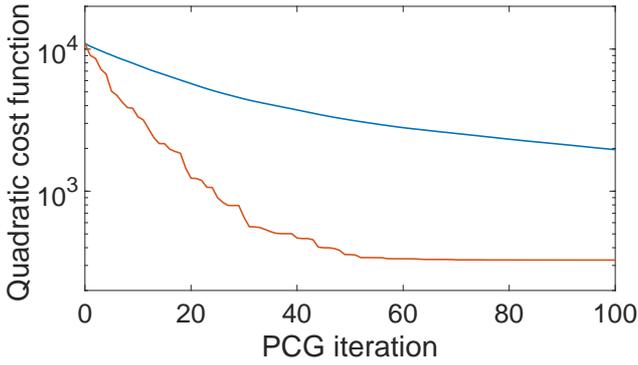}
  \caption{Case~\ref{case:few_obs}}
\end{subfigure}\\[1ex]
\caption{The values of the quadratic cost functions at every PCG iteration when using no preconditioner and preconditioning using exact $\bL^{-1}$. Values of $\sigma_o$ and the number of observations $p$ are varied in cases~\ref{case:reference}, \ref{case:large_sigmao}, \ref{case:few_obs}.}
\label{fig:state_forc_ps66-2_5_ft10_fx5_25}
\end{figure}

We have noticed that the effectiveness of the exact preconditioner $\bL^{-1}$ depends on how much of the system is observed and the interaction between the model and observation errors. There are observations at every 10th time step, ensuring that there are observations at the final time. We consider the following cases regarding the observation error variance $\sigma_o $ and the total number of observations $p$:
\begin{enumerate}
\item $\sigma_o = 1.5 \times 10^{-1}$, $p=300$ (observing $2\%$ of the system); \label{case:reference} 
\item $\sigma_o = 4.5 \times 10^{-1}$, $p=300$; \label{case:large_sigmao}
\item $\sigma_o = 1.5 \times 10^{-1}$, $p=60$ (observing $0.4\%$ of the system). \label{case:few_obs}
\end{enumerate}
In Figure~\ref{fig:state_forc_ps66-2_5_ft10_fx5_25}, we show that preconditioning using $\bL^{-1}$ is not useful in case~\ref{case:reference} but can be effective if the observation error variance is increased while keeping the same number of observations (case~\ref{case:large_sigmao}), or if the number of observations is reduced while $\sigma_o$ is unchanged (case~\ref{case:few_obs}).

Note that we compare the value of the quadratic cost function at every PCG iteration without taking into account the cost of the computation, which can be evaluated in terms of runtime or energy consumption and depends on how much parallelism can be achieved (e.g. \cite{CarsonStrakos2020}). If matrix-vector products with $\bL$ can be parallelised, then PCG iterations when solving the unpreconditioned system can be performed faster than with preconditioning. Then in terms of the runtime, preconditioning in case \ref{case:reference} can be even worse than indicated by comparing the quadratic cost function at every PCG iteration. In the same manner, preconditioning using exact $\bL^{-1}$ in cases \ref{case:large_sigmao} and \ref{case:few_obs} may not be as effective as displayed. In the following section, we test preconditioning using $\widetilde{\bL}^{-1}$ and $\widetilde{\bold{S}}$ in cases \ref{case:large_sigmao} and \ref{case:few_obs}.

\subsection{Preconditioning with randomised low rank approximation} 
\begin{figure}[h!]
\begin{subfigure}[b]{0.5\linewidth}
  \centering
 \includegraphics[width=\linewidth]{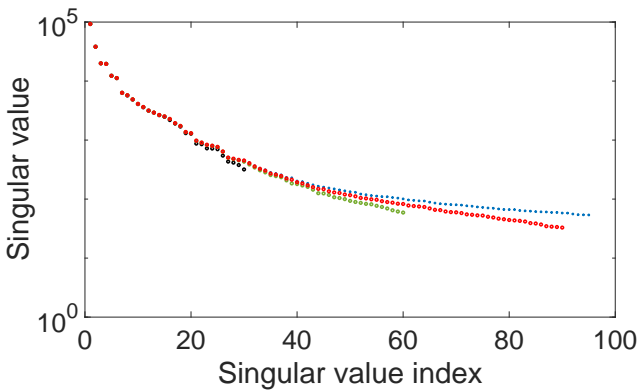}
\end{subfigure}
\begin{subfigure}[b]{0.5\linewidth}
  \centering
 \includegraphics[width=\linewidth]{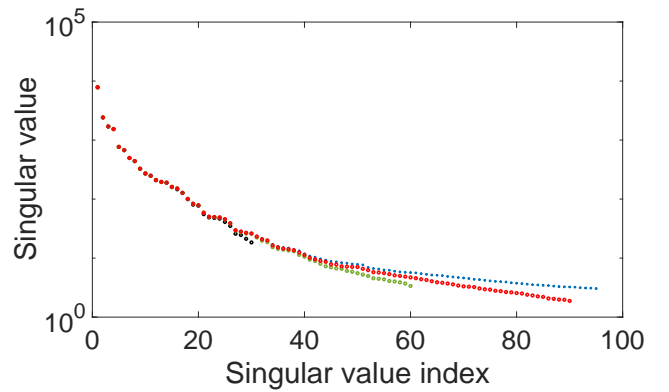}
\end{subfigure}\\[1ex]
\caption{Largest singular values of $\bold{P}$ (blue, left panel) and $\bold{W}$ (blue, right panel) and their approximations given by RSVD when using rank $k=30$ (black), $k=60$ (green) and $k=90$ (red). The largest singular values and their approximations coincide.}
\label{fig:rsvd_sing_val_k30_60_90_ps66-5_ft10_fx5}
\end{figure}

\begin{figure}[h!]
\begin{subfigure}[b]{0.5\linewidth}
  \centering
 \includegraphics[width=\linewidth]{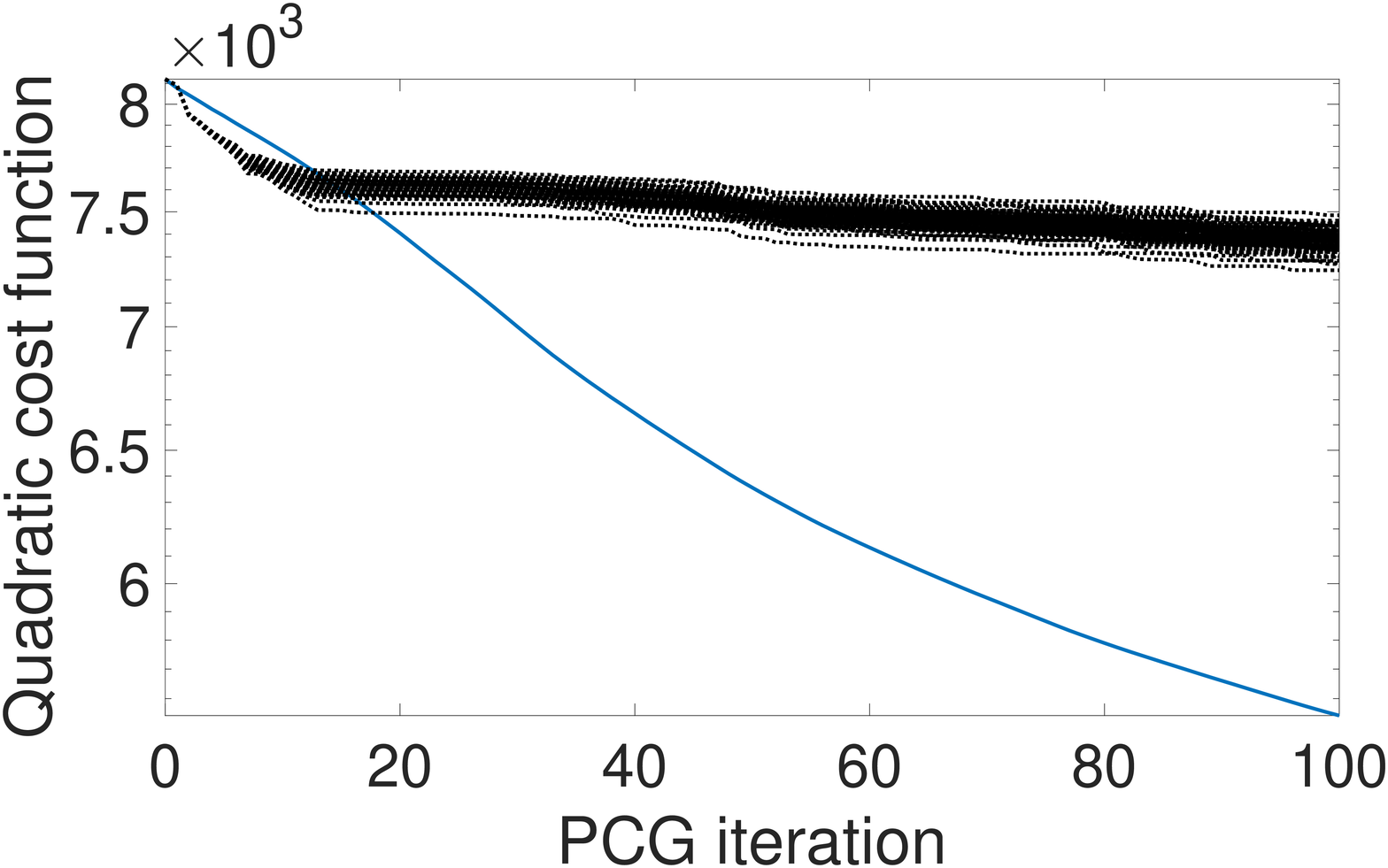}
  \caption{Case~\ref{case:large_sigmao}, $k=30$}
\end{subfigure}
\begin{subfigure}[b]{0.5\linewidth}
  \centering
 \includegraphics[width=\linewidth]{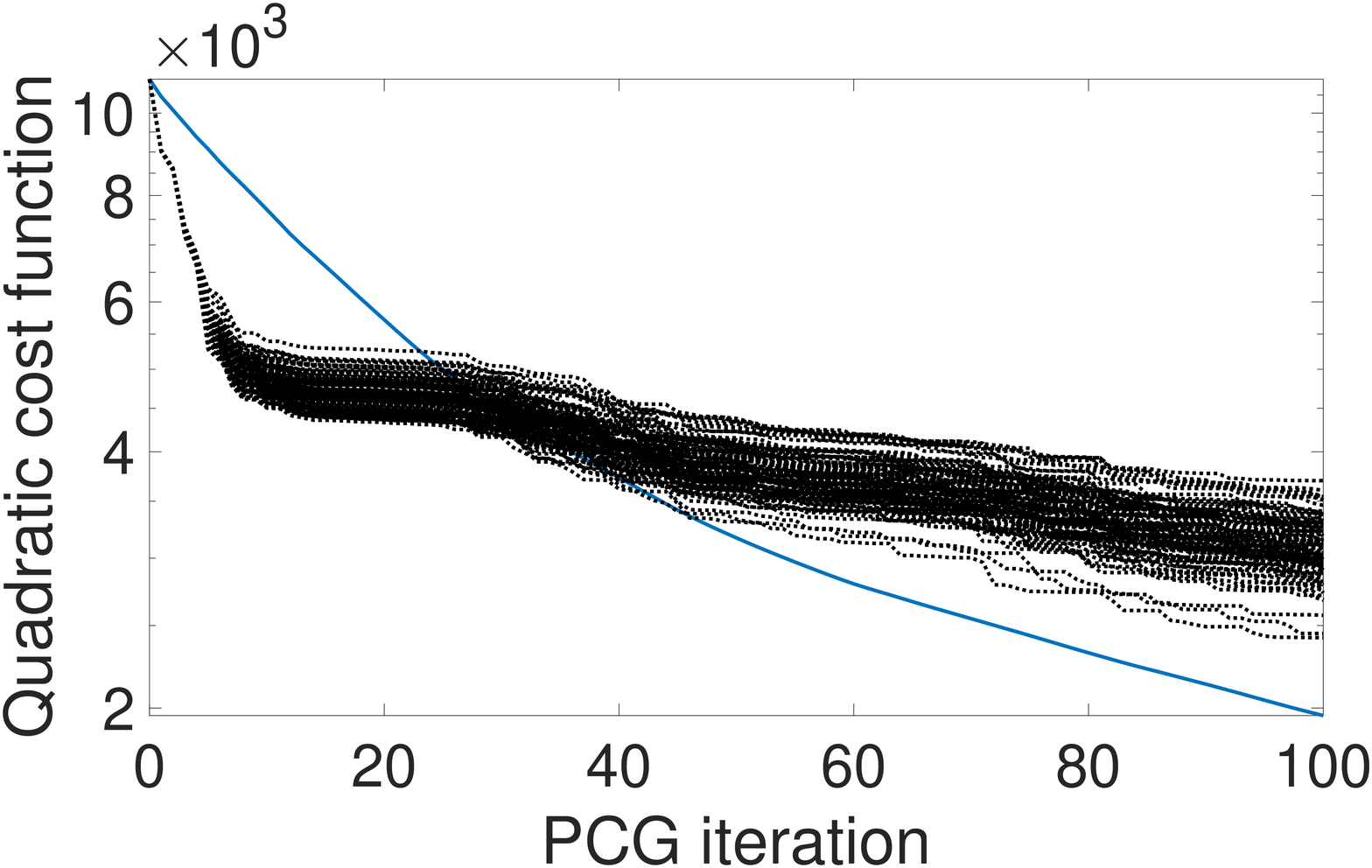}
  \caption{Case~\ref{case:few_obs}, $k=30$}
\end{subfigure}\\[1ex]
\begin{subfigure}[b]{0.5\linewidth}
  \centering
 \includegraphics[width=\linewidth]{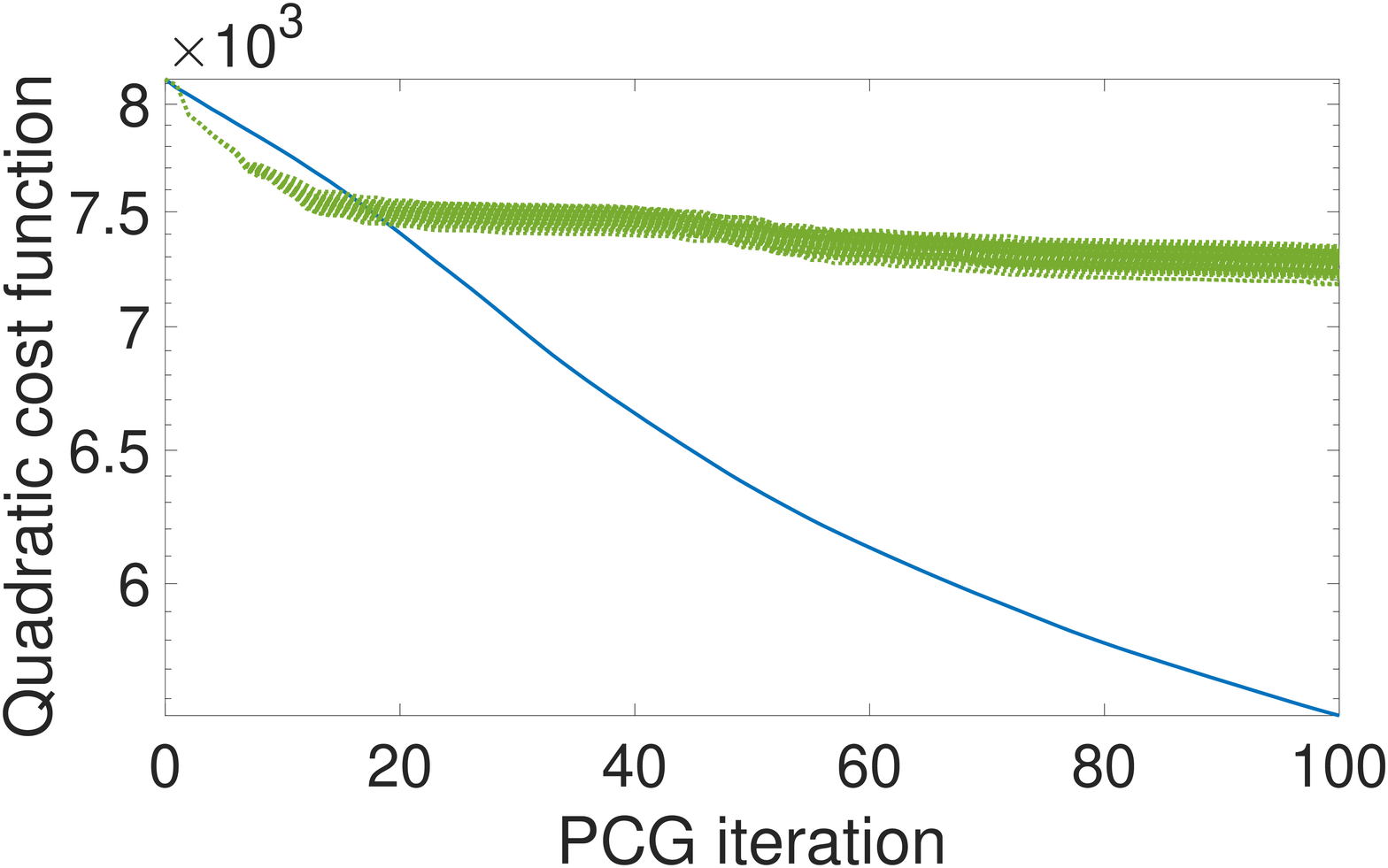}
  \caption{Case~\ref{case:large_sigmao}, $k=60$}
\end{subfigure}
\begin{subfigure}[b]{0.5\linewidth}
  \centering
 \includegraphics[width=\linewidth]{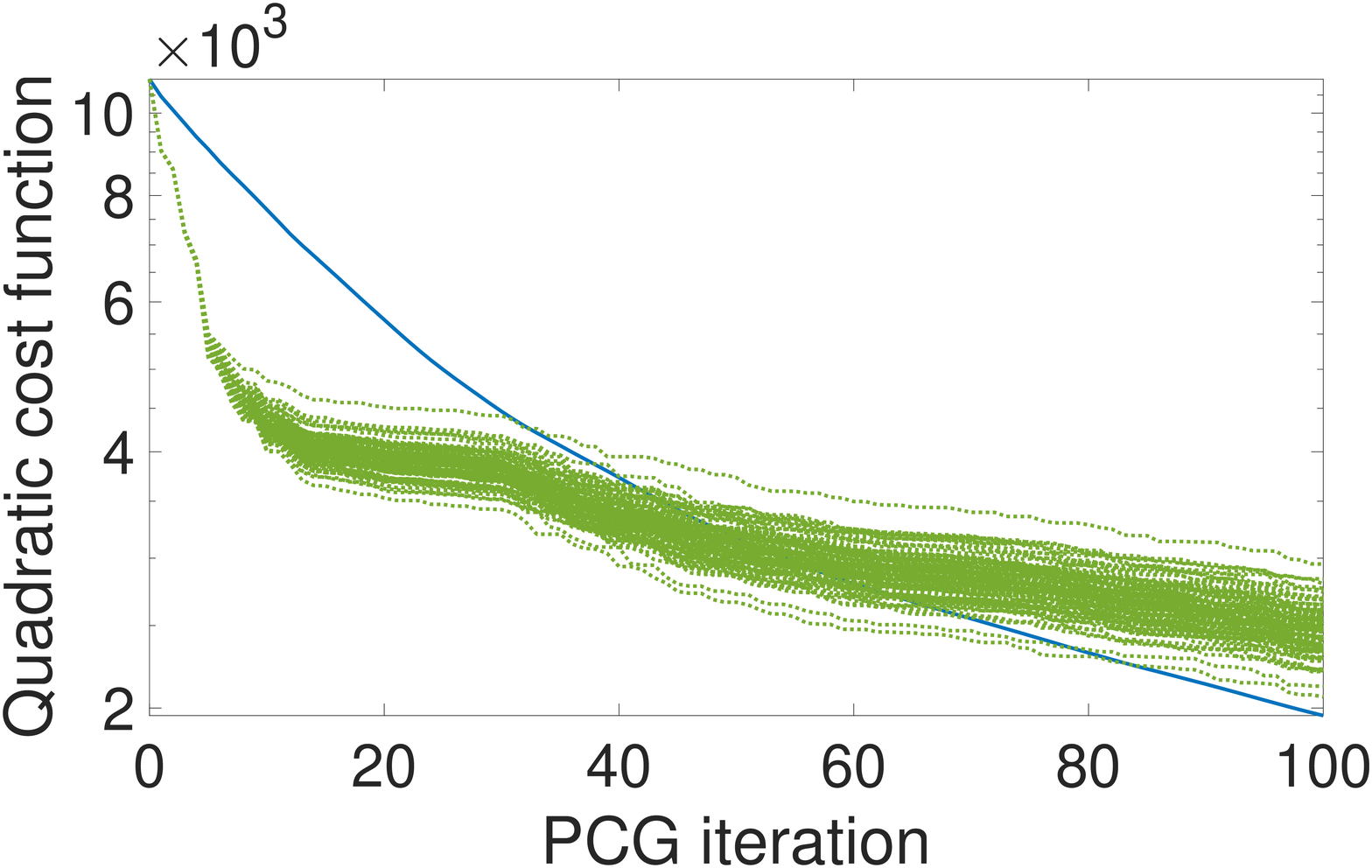}
  \caption{Case~\ref{case:few_obs}, $k=60$}
\end{subfigure}\\[1ex]
\begin{subfigure}[b]{0.5\linewidth}
  \centering
 \includegraphics[width=\linewidth]{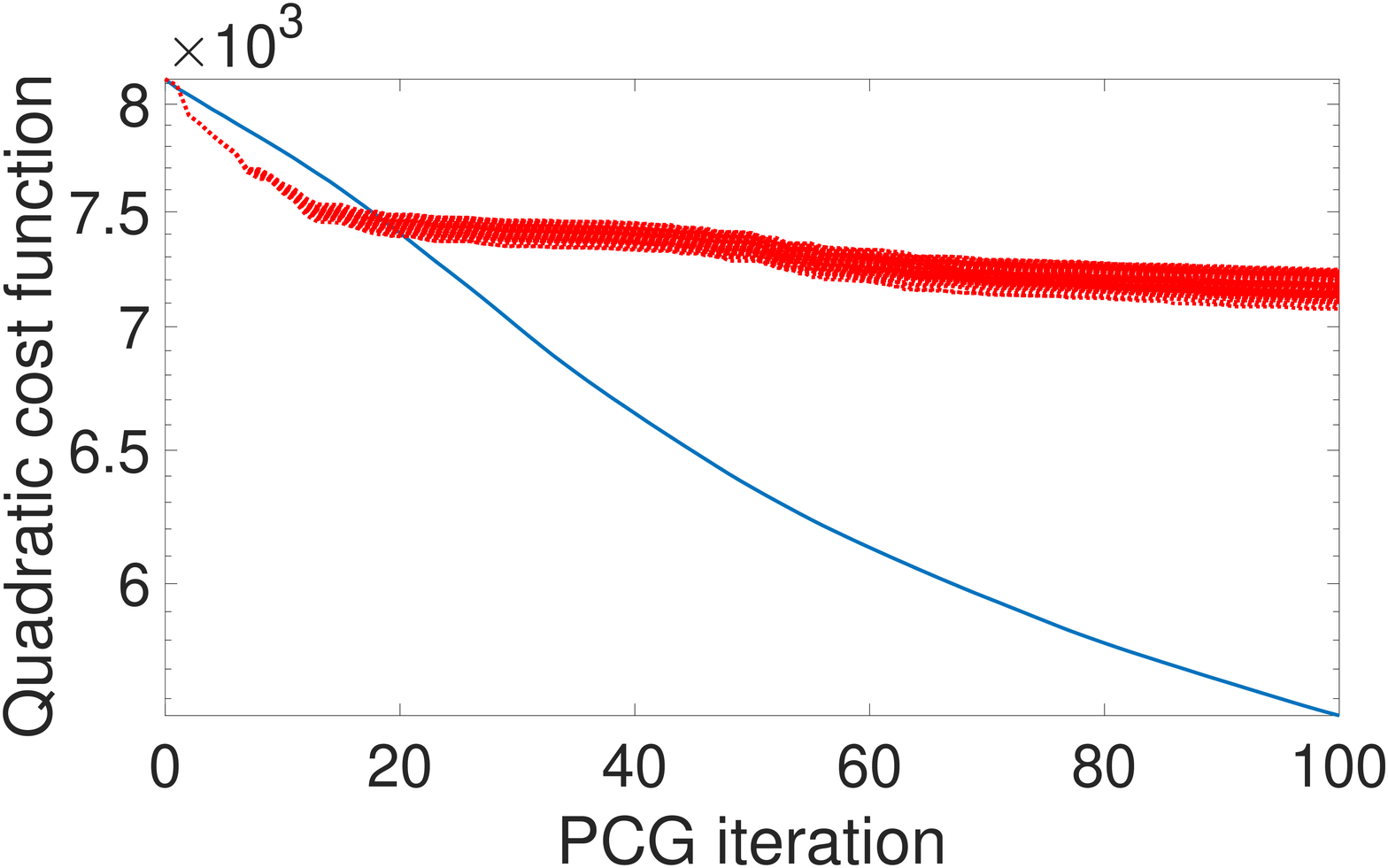}
  \caption{Case~\ref{case:large_sigmao}, $k=90$}
\end{subfigure}
\begin{subfigure}[b]{0.5\linewidth}
  \centering
 \includegraphics[width=\linewidth]{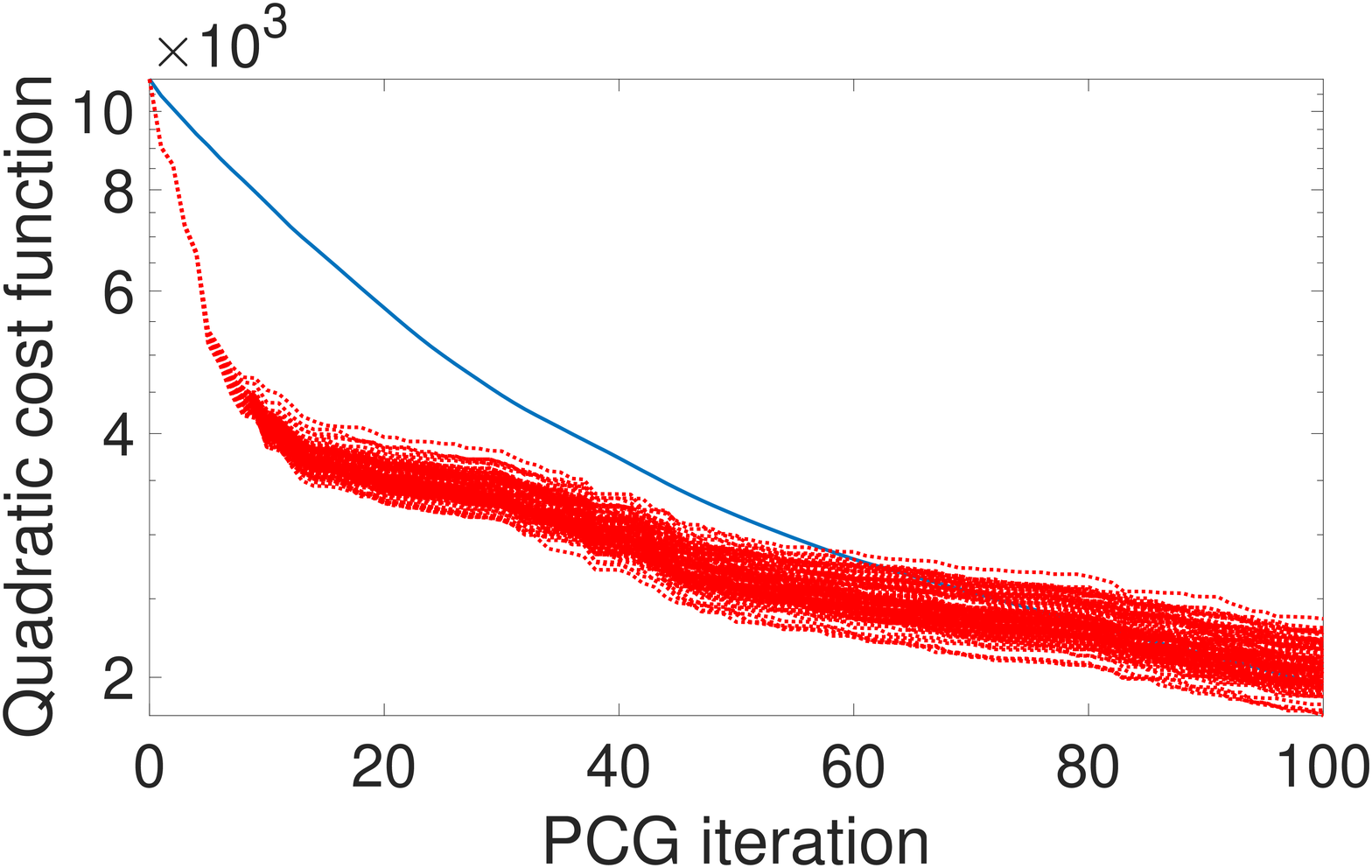}
  \caption{Case~\ref{case:few_obs}, $k=90$}
\end{subfigure}\\[1ex]
\caption{Values of the quadratic cost function at every PCG iteration when using no preconditioner (blue solid line) and preconditioning using $\widetilde{\bold{S}}$ (dotted lines) that are constructed using rank $k\in \{30,60,90\}$ approximation. One hundred realisations of the randomised preconditioner are shown. Values of $\sigma_o$ and the number of observations $p$ are varied in cases~\ref{case:large_sigmao} and \ref{case:few_obs}.}
\label{fig:preconditioned_all_cases_RV60_ps66-2_5_ft10_fx5_25}
\end{figure}

We generate $\widetilde{\bold{L}}^{-1}$ and $\widetilde{\bold{S}}$ by using rank $k \in \{30, 60, 90 \}$ approximations of $\bold{P}$ and $\bold{W}$ in \eqref{eq:Linv_I_plus_P} and \eqref{eq:S_Dsqrt_plus_W}, respectively. The oversampling parameter is set to $l=5$. We found that using $l=10$ or $l=15$ does not make a significant difference to the results (not shown). RSVD produces high quality approximations of the singular values of both $\bold{P}$ and $\bold{W}$. The largest singular values and their approximations are shown in Figure~\ref{fig:rsvd_sing_val_k30_60_90_ps66-5_ft10_fx5}, where the same random seed is used to generate the random matrix $\bold{G}$ for all $k$ values. Matrices $\bold{P}$ and $\bold{W}$ do not depend on whether case~\ref{case:large_sigmao} or \ref{case:few_obs} is considered, because the cases differ in the observation terms. In each case, we run the RSVD algorithm one hundred times with different Gaussian matrices $\bold{G}$ and solve the systems with the resulting preconditioners. The spread is illustrated in Figure~\ref{fig:preconditioned_all_cases_RV60_ps66-2_5_ft10_fx5_25} for $\widetilde{\bold{S}}$. In both cases, the variation in the values of the cost function is small during the early iterations. This shows that our results are not very sensitive to the choice of $\bold{G}$ and, in practice, it is only necessary to run the RSVD algorithm once. 

\begin{figure}[h!]
\begin{subfigure}[b]{0.5\linewidth}
  \centering
 \includegraphics[width=\linewidth]{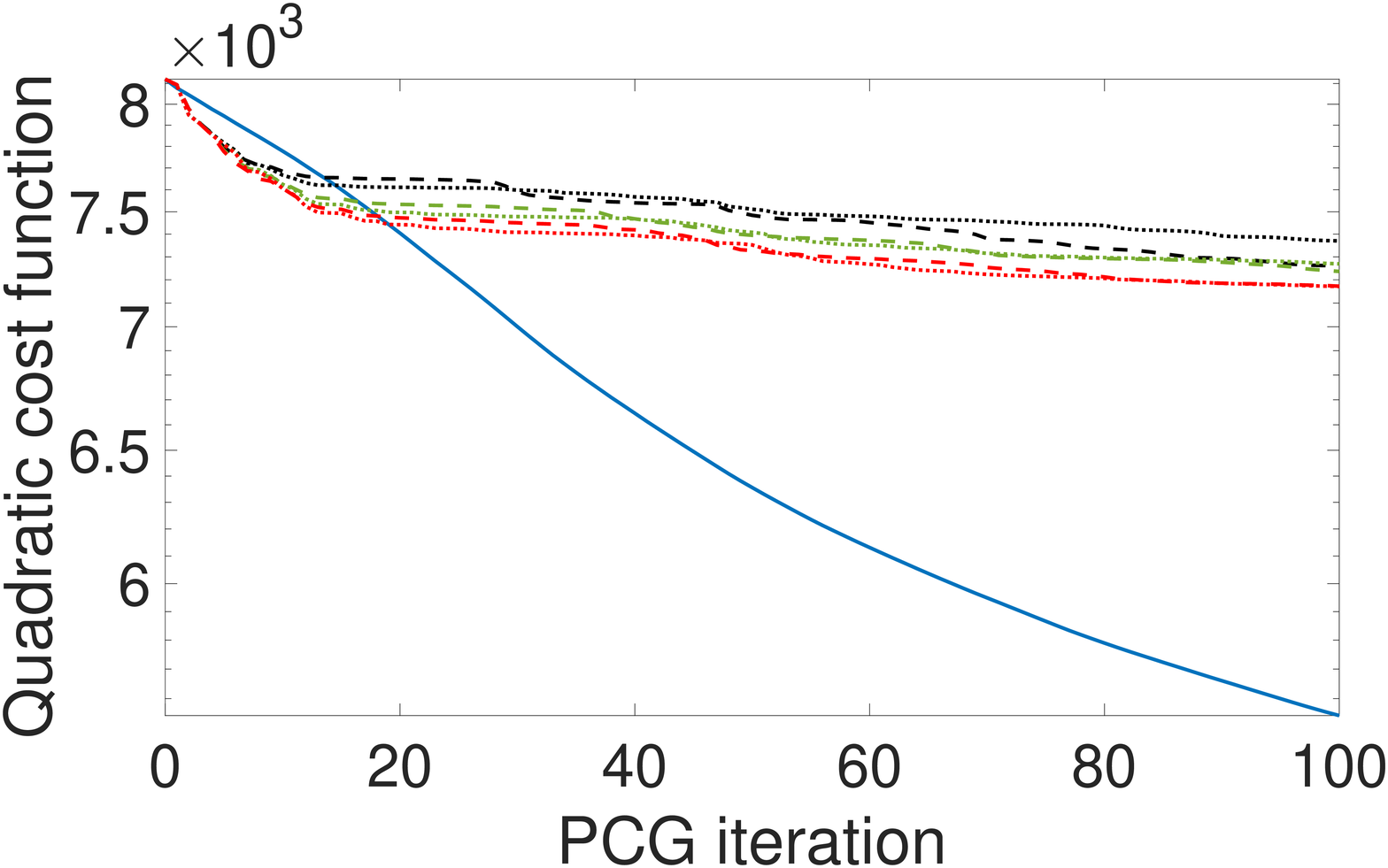}
  \caption{Case~\ref{case:large_sigmao}}
\end{subfigure}
\begin{subfigure}[b]{0.5\linewidth}
  \centering
 \includegraphics[width=\linewidth]{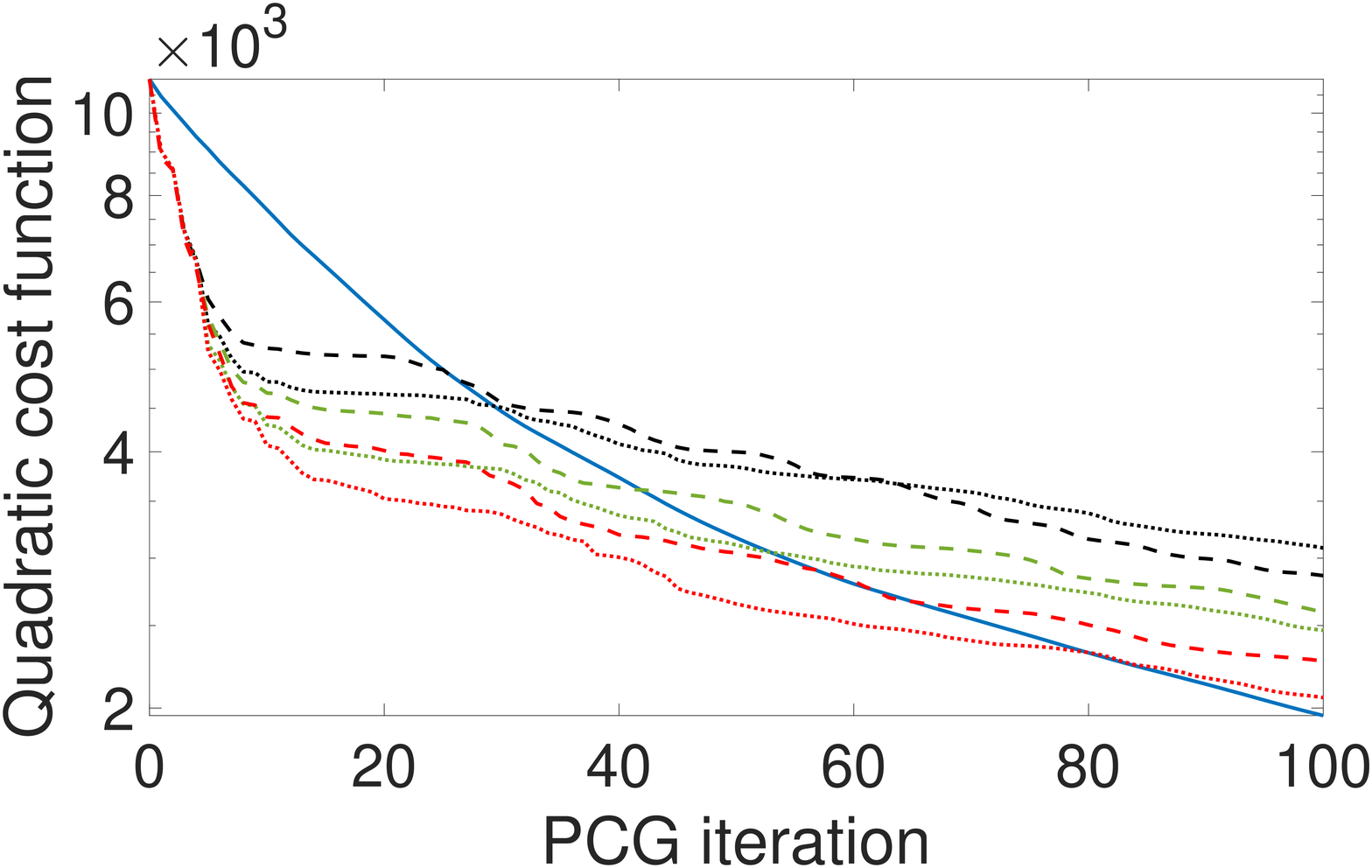}
  \caption{Case~\ref{case:few_obs}}
\end{subfigure}\\[1ex]
\caption{Mean values (over one hundred experiments) of the quadratic cost function at every PCG iteration when using no preconditioner (blue solid line) and preconditioning using $\widetilde{\bold{L}}^{-1}$ (dashed lines) and $\widetilde{\bold{S}}$ (dotted lines) that are constructed using rank $k=30$ (black), $k=60$ (green) and $k=90$ (red) approximation. Values of $\sigma_o$ and the number of observations $p$ are varied in cases~\ref{case:large_sigmao} and \ref{case:few_obs}.}
\label{fig:preconditioned_cases_ps66-2_5_ft10_fx5_25}
\end{figure}

The means of the quadratic cost function in cases~\ref{case:large_sigmao} and \ref{case:few_obs} are shown in Figure~\ref{fig:preconditioned_cases_ps66-2_5_ft10_fx5_25}. Higher rank approximations in both cases and using $\widetilde{\bold{S}}$ in case~\ref{case:few_obs} results in faster minimisation. Notice that in the first few iterations of PCG, preconditioning gives the same improvement regardless of the rank of approximation and whether $\widetilde{\bold{L}}^{-1}$ or $\widetilde{\bold{S}}$ is used. Preconditioning is more useful in case~\ref{case:few_obs}, which has fewer observations. The approximations used to generate $\widetilde{\bold{L}}^{-1}$ and $\widetilde{\bold{S}}$ are very low rank compared to the size of the system and there is a good improvement over the unpreconditioned case when the number of observations is low, especially in the beginning of the iterative process, which is the most relevant in practical settings. In the case with more observations (case~\ref{case:large_sigmao}), the randomised preconditioning is useful if a small number of PCG iterations is run. Since in an operational context we only run a small number of iterations, we are more likely to be in this regime. In cases~\ref{case:large_sigmao} and \ref{case:few_obs}, using exact $\bold{L}^{-1}$ results in a modest (case~\ref{case:large_sigmao}) and a rapid (case~\ref{case:few_obs}) decrease of the cost function in the first PCG iterations (Figure~\ref{fig:state_forc_ps66-2_5_ft10_fx5_25}). Our proposed preconditioners replicate such behaviour and if larger $k$ is used then the performance of exact $\bold{L}^{-1}$ is followed for more PCG iterations. In case~\ref{case:few_obs}, the quadratic cost function value is reduced by a factor of two after 5 PCG iterations when using exact $\bL^{-1}$ in the preconditioner, the same result is obtained after 8 ($k=30$) and 6 ($k=60$ and $k=90$) PCG iterations using $\widetilde{\bold{L}}^{-1}$, and 6 ($k=30$) and 5 ($k=60$ and $k=90$) PCG iterations using $\widetilde{\bold{S}}$. In case~\ref{case:large_sigmao}, the quadratic cost function is reduced only by a factor of 1.7 in 100 PCG iterations when preconditioning with the exact $\bL^{-1}$. When using our preconditioners the values of the quadratic cost function after 100 PCG iterations are larger than when using exact $\bL^{-1}$ or no preconditioning. This can be addressed by using a larger rank approximation, computational resources permitting.

\subsubsection{Large model error}
\begin{figure}[h!]
\begin{subfigure}[b]{0.5\linewidth}
  \centering
 \includegraphics[width=\linewidth]{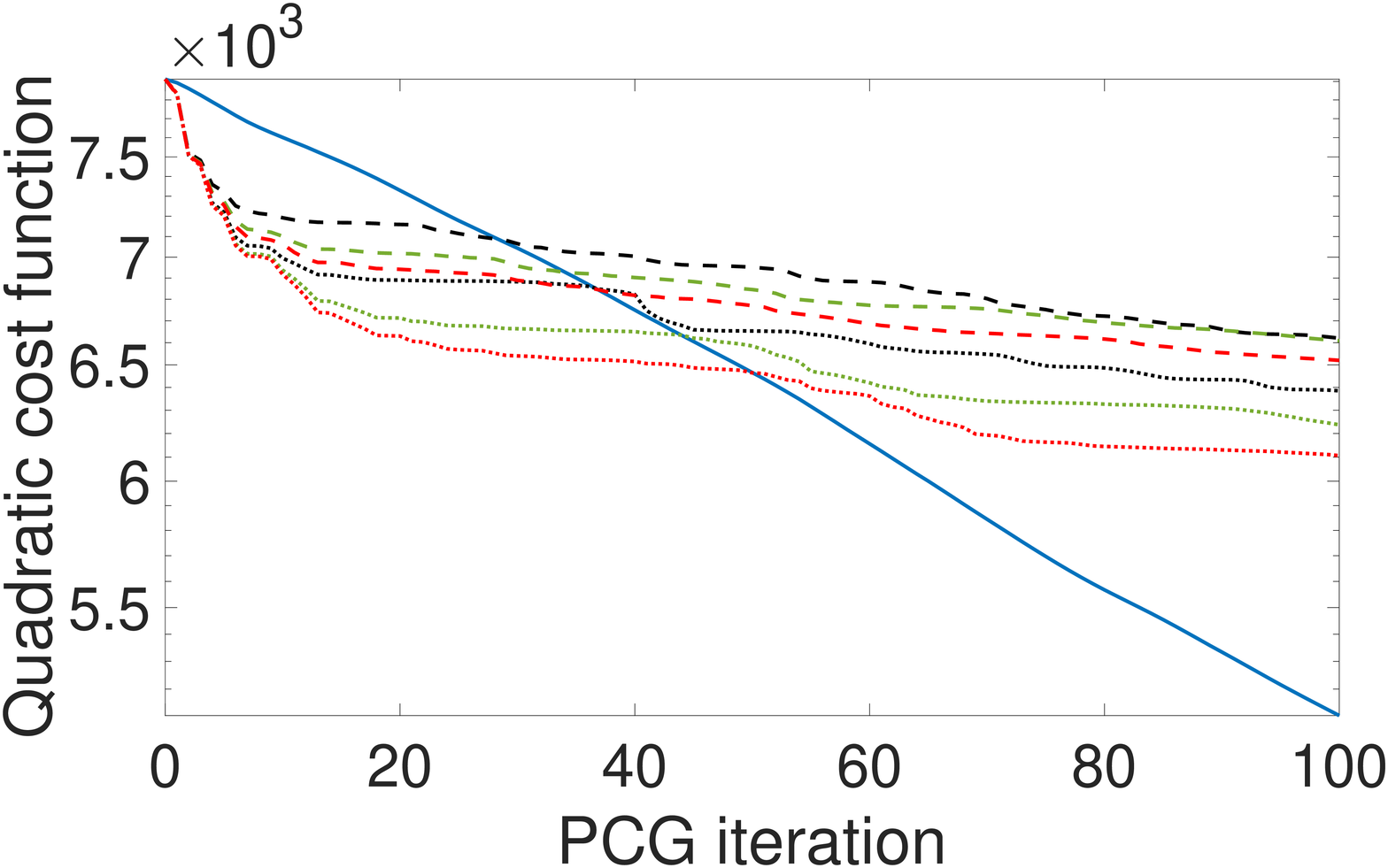}
  \caption{Case~\ref{case:large_sigmao}}
\end{subfigure}
\begin{subfigure}[b]{0.5\linewidth}
  \centering
 \includegraphics[width=\linewidth]{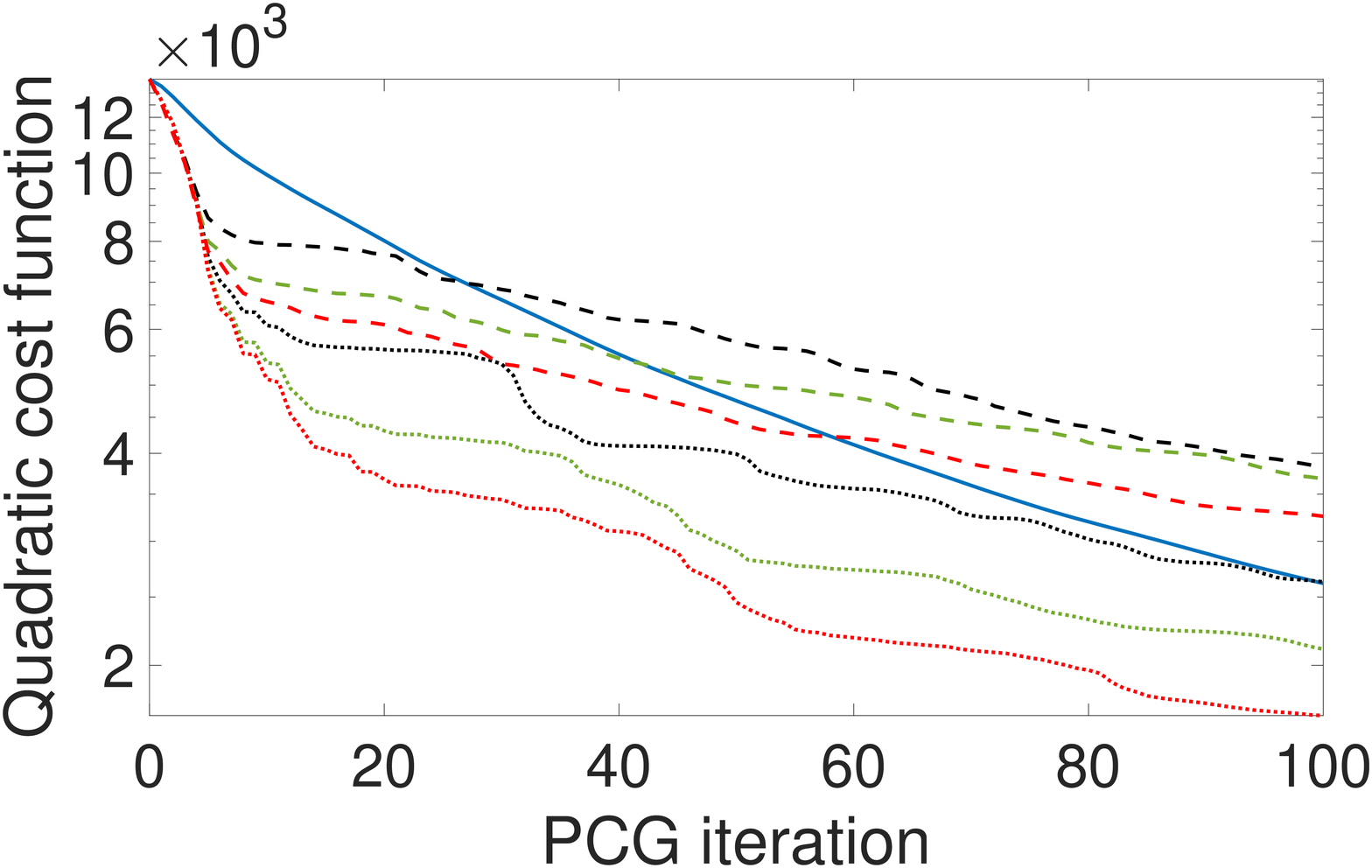}
  \caption{Case~\ref{case:few_obs}}
\end{subfigure}\\[1ex]
\caption{As in Figure~\ref{fig:preconditioned_cases_ps66-2_5_ft10_fx5_25}, but the model error covariance matrix is $\bQ_i = 0.1^2 \bC_q$ and $\bC_b$ has length scale $2\Delta X$.}
\label{fig:preconditioned_cases_ps37-2_5_ft10_fx5_25}
\end{figure}

We explore how the preconditioning using approximations of $\bL^{-1}$ and $\bL^{-1} \bold{D}^{1/2}$ compare when the model error is large. The numerical experiments are performed using the same setup as before, but now we set $\bold{Q}_i = 0.1^2 \bold{C}_q$, where $\bC_b$ has length scale $2\Delta X$. The means over one hundred runs are presented in Figure~\ref{fig:preconditioned_cases_ps37-2_5_ft10_fx5_25}. There is a clear separation between the minimisation using $\widetilde{\bold{L}}^{-1}$ and $\widetilde{\bold{S}}$ in the preconditioner after the first few PCG iterations; with the latter resulting in faster minimisation. Notice that the preconditioning using both approximations remains useful for more PCG iterations than in the setup with a smaller model error. This can be expected because the increase of length scales of $\bQ_i$ has a detrimental effect on the conditioning of the unpreconditioned Hessian (see, e.g. Chapter 6 of \cite{AES_thesis}) and hence preconditioning can be more efficient.

\section{Conclusions}\label{sec:conclusions}
We have considered preconditioning for the state formulation of incremental weak constraint 4D-Var, which closely follows the control variable transform (first level preconditioning) strategy for the strong constraint formulation. We have shown that such preconditioning may not be useful even when using the exact $\bL^{-1}$, which also makes the matrix-vector products with the Hessian sequential in the time dimension. In the cases where such preconditioning is useful, a good preconditioner can be obtained by using randomised singular value decompositions to approximate $\bL^{-1}$ or $\bL^{-1} \bold{D}^{1/2}$. These preconditioners are cheap to compute and apply and do allow for parallelization in the time dimension. They can improve the solution of the exact inner loop problem, resulting in a higher reduction of the quadratic cost function in the same number of iterations compared to using no preconditioning or obtain the same quadratic cost function value in fewer iterations. The effect of the accuracy of the inner loop solution on the analysis has been studied by e.g. \cite{Lawless2006}.

Our results call for caution when designing preconditioning approaches that focus on approximating $\bL^{-1}$, especially when the number of observations is high. In practical NWP settings, around $1\%$ of the system is observed, hence approximating $\bL^{-1}$ may be useful. Using randomised approximations of $\bL^{-1} $ or $\bL^{-1} \bold{D}^{1/2}$ should be tested using large and more realistic systems, where meaningful evaluations of the runtime and energy consumption can be obtained. A more detailed investigation on when preconditioning with $\bL^{-1}$ gives good results would also be useful.

\section*{Acknowledgements}
We thank Dr. Adam El-Said for his code for the weak constraint 4D-Var assimilation system. We are also grateful for the helpful comments by two anonymous reviewers.

\section*{Conflict of interest}
The authors declare no conflict of interest.

\section*{Funding information}
UK Engineering and Physical Sciences Research Council, Grant/Award Number: EP/L016613/1; European Research Council CUNDA project, Grant/Award Number: 694509; NERC National Centre for Earth Observation.

\bibliographystyle{apalike}
\bibliography{../../../../PhD_biblio}

\begin{thebibliography}{24}
\providecommand{\natexlab}[1]{#1}
\providecommand{\url}[1]{\texttt{#1}}
\expandafter\ifx\csname urlstyle\endcsname\relax
  \providecommand{\doi}[1]{doi: #1}\else
  \providecommand{\doi}{doi: \begingroup \urlstyle{rm}\Url}\fi

\bibitem[Bousserez et~al.(2020)Bousserez, Guerrette, and Henze]{Bousserez2020}
N.~Bousserez, J.~J. Guerrette, and D.~K. Henze.
\newblock Enhanced parallelization of the incremental 4{D}-{V}ar data
  assimilation algorithm using the {R}andomized {I}ncremental {O}ptimal
  {T}echnique.
\newblock \emph{Quarterly Journal of the Royal Meteorological Society},
  146\penalty0 (728):\penalty0 1351 -- 1371, 2020.

\bibitem[Butcher(1987)]{Butcher87}
J.~C. Butcher.
\newblock \emph{The Numerical Analysis of Ordinary Differential Equations:
  Runge-Kutta and General Linear Methods}.
\newblock Wiley-Interscience, 1987.

\bibitem[Carson and Strako\v{s}(2020)]{CarsonStrakos2020}
E.~Carson and Z.~Strako\v{s}.
\newblock On the cost of iterative computations.
\newblock \emph{Philosophical Transactions of the Royal Society A:
  Mathematical, Physical and Engineering Sciences}, 378\penalty0
  (2166):\penalty0 20190050, 2020.

\bibitem[Daley(1993)]{Daley91}
R.~Daley.
\newblock \emph{Atmospheric Data Analysis}, volume~2.
\newblock Cambridge University Press, 1993.

\bibitem[Dau\v{z}ickait\.{e} et~al.(2021)Dau\v{z}ickait\.{e}, Lawless, Scott,
  and van Leeuwen]{Dauzickaite2021}
I.~Dau\v{z}ickait\.{e}, A.~S. Lawless, J.~A. Scott, and P.~J. van Leeuwen.
\newblock Randomised preconditioning for the forcing formulation of weak
  constraint 4{D}-{V}ar.
\newblock 2021.
\newblock Submitted. Available on \url{https://arxiv.org/abs/2101.07249}.

\bibitem[El-Said(2015)]{AES_thesis}
A.~El-Said.
\newblock \emph{Conditioning of the weak-constraint variational data
  assimilation problem for numerical weather prediction}.
\newblock PhD thesis, Department of Mathematics and Statistics, University of
  Reading, 2015.

\bibitem[Fisher and G{\"u}rol(2017)]{Fisher17}
M.~Fisher and S.~G{\"u}rol.
\newblock Parallelisation in the time dimension of four-dimensional variational
  data assimilation.
\newblock \emph{Quarterly Journal of the Royal Meteorological Society},
  143\penalty0 (703):\penalty0 1136 -- 1147, 2017.

\bibitem[Gratton et~al.(2007)Gratton, Lawless, and Nichols]{Gratton2007}
S.~Gratton, A.~S. Lawless, and N.~K. Nichols.
\newblock Approximate {G}auss-{N}ewton methods for nonlinear least squares
  problems.
\newblock \emph{SIAM Journal on Optimization}, 18\penalty0 (1):\penalty0 106 --
  132, 2007.

\bibitem[Gratton et~al.(2018{\natexlab{a}})Gratton, G{\"u}rol, Simon, and
  Toint]{Gratton18}
S.~Gratton, S.~G{\"u}rol, E.~Simon, and Ph.~L. Toint.
\newblock A note on preconditioning weighted linear least-squares, with
  consequences for weakly constrained variational data assimilation.
\newblock \emph{Quarterly Journal of the Royal Meteorological Society},
  144\penalty0 (712):\penalty0 934 -- 940, 2018{\natexlab{a}}.

\bibitem[Gratton et~al.(2018{\natexlab{b}})Gratton, G{\"{u}}rol, Simon, and
  Toint]{Gratton2018}
S.~Gratton, S.~G{\"{u}}rol, E.~Simon, and Ph.~L. Toint.
\newblock Guaranteeing the convergence of the saddle formulation for
  weakly-constrained 4{D}-{V}ar data assimilation.
\newblock \emph{Quarterly Journal of the Royal Meteorological Society},
  144\penalty0 (717):\penalty0 2592 -- 2602, 2018{\natexlab{b}}.

\bibitem[Gu(2015)]{Gu2015}
M.~Gu.
\newblock Subspace iteration randomization and singular value problems.
\newblock \emph{SIAM Journal on Scientific Computing}, 37\penalty0
  (3):\penalty0 A1139 -- A1173, 2015.

\bibitem[Halko et~al.(2011)Halko, Martinsson, and Tropp]{Halko11}
N.~Halko, P.~Martinsson, and J.~Tropp.
\newblock Finding structure with randomness: Probabilistic algorithms for
  constructing approximate matrix decompositions.
\newblock \emph{SIAM Review}, 53\penalty0 (2):\penalty0 217 -- 288, 2011.

\bibitem[Johnson et~al.(2005)Johnson, Hoskins, and Nichols]{Johnson2005}
C.~Johnson, B.~J. Hoskins, and N.~K. Nichols.
\newblock A singular vector perspective of 4{D}-{V}ar: {F}iltering and
  interpolation.
\newblock \emph{Quarterly Journal of the Royal Meteorological Society},
  131\penalty0 (605):\penalty0 1 -- 19, 2005.

\bibitem[Lawless(2013)]{Lawless2013}
A.~S. Lawless.
\newblock Variational data assimilation for very large environmental problems.
\newblock In M.J.P Cullen, M.~A. Freitag, S.~Kindermann, and R.~Scheichl,
  editors, \emph{Large Scale Inverse Problems: Computational Methods and
  Applications in the Earth Sciences}, Radon Series on Computational and
  Applied Mathematics 13, pages 55 -- 90. De Gruyter, 2013.

\bibitem[Lawless and Nichols(2006)]{Lawless2006}
A.~S. Lawless and N.~K. Nichols.
\newblock Inner-loop stopping criteria for incremental four-dimensional
  variational data assimilation.
\newblock \emph{Monthly Weather Review}, 134\penalty0 (11):\penalty0 3425 --
  3435, 2006.

\bibitem[Lawless et~al.(2008)Lawless, Nichols, Boess, and
  Bunse-Gerstner]{Lawless2008}
A.~S. Lawless, N.~K. Nichols, C.~Boess, and A.~Bunse-Gerstner.
\newblock Using model reduction methods within incremental 4d-var.
\newblock \emph{Monthly Weather Review}, 136:\penalty0 1511 -- 1522, 2008.

\bibitem[Liesen and Strako\v{s}(2013)]{Liesen_Strakos}
J.~Liesen and Z.~Strako\v{s}.
\newblock \emph{Krylov subspace methods: principles and analysis. Numerical
  Mathematics and Scientific Computation}.
\newblock Oxford University Press, 2013.

\bibitem[Lorenc et~al.(2000)Lorenc, Ballard, Bell, Ingleby, Andrews, Barker,
  Bray, Clayton, Dalby, Li, Payne, and Saunders]{Lorenc2000}
A.~C. Lorenc, S.~P. Ballard, R.~S. Bell, N.~B. Ingleby, P.~L.~F. Andrews, D.~M.
  Barker, J.~R. Bray, A.~M. Clayton, T.~Dalby, D.~Li, T.~J. Payne, and F.~W.
  Saunders.
\newblock The {M}et {O}ffice global three-dimensional variational data
  assimilation scheme.
\newblock \emph{Quarterly Journal of the Royal Meteorological Society},
  126\penalty0 (570):\penalty0 2991 -- 3012, 2000.

\bibitem[Lorenz(1996)]{Lorenz96}
E.~Lorenz.
\newblock Predictability - a problem partly solved.
\newblock In \emph{Proceedings of the Seminar on Predictability}, volume~1,
  pages 1 -- 18. European Centre for Medium Range Weather Forecasts, 1996.

\bibitem[Martinsson and Tropp(2020)]{Martinsson20}
P.~G. Martinsson and J.~A. Tropp.
\newblock Randomized numerical linear algebra: Foundations and algorithms.
\newblock \emph{Acta Numerica}, 29:\penalty0 403 -- 572, 2020.

\bibitem[Nocedal and Wright(2006)]{Nocedal06}
J.~Nocedal and S.~J. Wright.
\newblock \emph{Numerical Optimization}.
\newblock World Scientific, 2nd edition, 2006.

\bibitem[Rawlins et~al.(2007)Rawlins, Ballard, Bovis, Clayton, Li, Inverarity,
  Lorenc, and Payne]{Rawlins2007}
F.~Rawlins, S.~P. Ballard, K.~J. Bovis, A.~M. Clayton, D.~Li, G.~W. Inverarity,
  A.~C. Lorenc, and T.~J. Payne.
\newblock The {M}et {O}ffice global four-dimensional variational data
  assimilation scheme.
\newblock \emph{Quarterly Journal of the Royal Meteorological Society},
  133\penalty0 (623):\penalty0 347 -- 362, 2007.

\bibitem[Saad(2003)]{SaadBook}
Y.~Saad.
\newblock \emph{Iterative Methods for Sparse Linear Systems}.
\newblock {S}{I}{A}{M}, 2nd edition, 2003.

\bibitem[Tr{\'{e}}molet(2006)]{Tremolet06}
Y.~Tr{\'{e}}molet.
\newblock Accounting for an imperfect model in 4{D}-{V}ar.
\newblock \emph{Quarterly Journal of the Royal Meteorological Society},
  132\penalty0 (621):\penalty0 2483 -- 2504, 2006.

\end{thebibliography}

\end{document}